\documentstyle[11pt]{article}
\newtheorem{thm}{\sc Theorem}[section]
\newtheorem{prop}[thm]{\sc Proposition}
\newtheorem{lem}[thm]{\sc Lemma}

\newcommand{\sm}{\subset_{\!\!{ \rm m}}}
\newtheorem{rem}{\sc Remark}[section]
\newcommand{\m}{ {\rm mt}}
\newcommand{\mt}{minimally transitive }
\newcommand{\ka}{K_{G:A}}

\newcommand{\Om}{\Omega}

\newcommand{\pf}{\smallskip{\it Proof:} \,\,\,}
\newcommand{\dne}{\hfill $\Box$}
\begin{document}
\title{\sc  On Solvable Minimally Transitive  \\Permutation
Groups}

\author{ Francesca Dalla Volta\thanks{Work partially supported by M.I.U.R. and London Mathematical Society} \\  
{\em Dipartimento di Matematica e Applicazioni, Universit\`{a} Milano-Bicocca}\\
{\em 20125 Milano, Italy}\\
{\small\tt Francesca.DallaVolta@unimib.it}\\\\
and Johannes Siemons\\
{\em School of Mathematics, University of East Anglia,}\\
{\em Norwich, NR4 7TJ, United Kingdom.}\\{\small\tt
j.siemons@uea.ac.uk}}
\date{}
%%%%
%%%%\scriptsize {Version of 20 September 2006, amended for referees comments on 30th January 2007printed \today
%%%%
\maketitle

\begin{abstract}
\noindent We investigate properties of finite transitive permutation groups
$(G,\, \Om)$ in which all proper subgroups of $G$ act
intransitively on $\Om.$ In particular,  we are interested in reduction theorems for
minimally transitive
representations of solvable groups.
\end{abstract}

\section{\sc  Introduction}

A finite permutation group $(G,\, \Om)$ is {\it minimally
transitive } if $G$ is transitive on $\Om$ while all proper
subgroups of $G$ are intransitive on $\Om.$ Evidently, any
transitive permutation group contains minimally transitive
subgroups acting on the same set and so this concept occurs
naturally in reduction arguments. Closely related is the notion of
minimally irreducible linear groups, namely those linear groups
which act irreducibly on a vector space  $V$ while all  proper
subgroups  leave some proper subspace of $V$  invariant.

Solvable minimally transitive groups were first considered by
Suprunenko~\cite{Sup1} and  Kopylova~\cite{Kopy1} who studied the
groups of degree $pq$ with $p$ and $q$ primes. More recently
Lucchini~\cite{Lucch1} studied minimal generating sets in
minimally transitive groups, in connection with the asymptotic
properties of permutation groups considered in
Pyber~\cite{Pyber1}. In Ngo~\cite{Dak} non-regular metabelian
minimally transitive groups are investigated, and
Miller-Praeger~\cite{Praeg1} mention such groups in the
context of vertex transitive graphs which are not Cayley graphs. A
list of minimally transitive groups up to degree $30$ is available
in Hulpke~\cite{Hulp1}, see also Conway, Hulpke and
McKay~\cite{Hulp2}.

In this paper we consider solvable groups. Here in particular it happens
frequently that a group action  is not faithful. Therefore we study more generally  arbitrary  {\it minimally transitive representations} which may or may not be faithful.  This language requires technical detail which could detract from the main matter; wherever possible we therefore try to stay close to the language of permutation groups which may appear more natural. 
  
Any transitive permutation group contains  minimally
transitive subgroups and therefore it would be unreasonable to expect full classifications in general. However, under suitable restrictions some general results can be expected. For instance, for  nilpotent groups there is a simple description of all their minimally transitive representations. 

For a faithful action the primes dividing the order of a solvable group must divide the degree, see Theorem~\ref{1012}.  In Sections 2 and 3  we prove some reduction theorems for subgroups and factor groups. In particular, a construction is given to reduce a general minimally transitive action to the case where the degree contains only two primes. A good result is also available for actions of square-free degree, extending the work of Suprunenko and  Kopylova. 
 
\section{\sc Minimally Transitive Groups}

Let $G\subseteq {\rm Sym\,}\Om$ be a transitive permutation  group
on a finite set  $\Om.$ Then $G$ is {\it minimally transitive} on
$\Om$ if every proper subgroup of $G$ is intransitive on $\Omega.$
In the following we consider
more generally an abstract finite group $G$  together with all its
transitive actions, faithful or not. Thus if $A\subset G$ is a
subgroup of $G$ let $G$ act  on the cosets $G\!:\!A$ of
$A$ in $G.$ The kernel of this action is the core
$\ka:=\bigcap_{g\in G}\,A^{g}$ of $A$ in $G.$

Thus $G$ acts {\it minimally transitively } on $G\!:\!A$ if and
only if every subgroup $H$ with $\ka\subseteq H \subset G$ acts
intransitively on $G\!:\!A.$ It will be convenient to call such a
subgroup $A$ an {\it mt-stabilizer} in $G$; we denote this as $A\sm
G.$ Therefore $A\sm G$ if and only if the following holds:
 Whenever $H\subseteq G$ is transitive on $G\!:\!A$ then $H\ka=G.$
 Evidently, if $A$ is an arbitrary subgroup of $G$ then $G/\ka$ always acts faithfully on $G\!:\!A$
 and hence is a permutation group on $G\!:\!A.$ This permutation group then is minimally transitive
 if and only if $A$ is an \m-stabilizer. 
 For instance, if $A=1$ then $G$ is  regular on $G\!:\!A$ and so $1\sm
G.$ For another example suppose that $|G|=pq$ with
distinct primes and Sylow subgroups $A\unlhd G$ and $B\,\,\,\,
/\!\!\!\!\!\!\unlhd G.$ Then $\ka=A$ and $A\sm G$ while 
$K_{G:B}=1$ and $B\not\sm G.$

 \vspace{5mm}
\subsection{\sc Preliminaries}

We begin by listing  general properties of groups with  \mt action.
For the remainder let $G$ be a finite group and let $A$ be a
subgroup of $G.$ The property of being an \m-stabilizer in $G$ is
quite special as it relates to the subgroup as well as its embedding in $G.$ Let
${\cal L}(G)$ denote the lattice of all subgroups of $G.$ We will
begin  by describing some  general properties of groups in ${\cal
L}(G)$ which are \m-stabilizers in  $G.$ The next lemma  is technical
but essential; the first part  we  use later on without further
mention.

\begin{lem}
\label{0009}
(i) \,\,\, Let $A\subseteq G.$ Then $A\sm G$ if and only if
$AH=G$ for a subgroup $H\subseteq G$ implies that $H\ka=G.$\newline (ii) \,\,
Let $A \sm G$ and let $B\subseteq A.$ Then  \,\,(a):\, $B\sm G$ or
\,\,(b):\, $K_{G:B}\neq \ka,$ $B\ka\sm G$ and there exists  a
subgroup $H\subseteq G$ with $HK_{G:B}\neq H\ka=G.$ In particular,
if
 $A\sm G$ and $\ka\subseteq B\subseteq A$  then $B\sm
G.$
\end{lem}

\pf (i) Suppose that $A\subseteq G$ and that also $H\subseteq G.$
Then $H$ is transitive on $G\!:\!A$ if and only if $AH=G.$
Therefore by definition, if $A\sm G$ and if $AH=G$ then $H\ka=G.$
Conversely, if $AH=G$ implies that $H\ka=G$ then $H$ being
transitive on $G\!:\!A$ means that $AH=G$ and so $H\ka=G.$ Hence
$A\sm G.$\newline
 (ii) Assume that $A \sm G$ and $B\subseteq A.$ If $B\not\sm G$ then there exists some $H$
such that $G=BHK_{G:B}$ but $G\neq HK_{G:B}.$ As $K_{G:B}\unlhd
\ka$ we have  $G=AH\ka.$ If $H\ka\neq G$ then $A\not\sm G.$ Next
we compute the core $\bar K$ of $B\ka$ in $G.$ Evidently,
$\ka\subseteq \bar K\subseteq B\ka\subseteq A$ so that $\bar
K=\ka.$ If $B\ka\not\sm G$ then there exists a subgroup $\bar H$
such that $(B\ka )\bar H=G$ but $\bar H\bar K\neq G.$ But then
$A(\ka \bar H)=G$ and $\bar H\ka\neq G,$ a contradiction. If
$\ka\subseteq B\subseteq A$ then $K_{G:B}=\ka$ and hence the
second alternative  can not happen. \dne

\bigskip
 When dealing with the set of all \m-stabilizers in
the subgroup lattice of  $G$ the following is a useful notion: If
$({\cal L},\,\leq)$ is any partially ordered set  then a subset
${\cal M}$ of ${\cal L}$ is an {\it order ideal\,} in ${\cal L}$
if $X\in {\cal M}$ and $Y\leq X$ with $Y\in{\cal L}$ implies that
$Y\in {\cal M}.$
\begin{rem}
\label{5001} From the last part of the lemma we 
deduce that the core-free {\rm mt}-stabilizers in $G$ form an
order ideal in the subgroup lattice ${\cal L}(G).$
\end{rem}

\medskip
It is therefore  often sufficient to know the 'top'
\m-stabilizers, that is those which are maximal subject to being
an \m-stabilizer. For instance, if $G$ is simple then the
\m-stabilizers form an order ideal and this is described
completely by its top elements. We may also  ask when such top
elements are  maximal subgroups of $G.$ Evidently,  $A$ is maximal
in  $G$ precisely when $G$ acts primitively on $G\!:\!A.$ More
generally, $G$ acts  {\it quasi-primitively} on  $G\!:\!A$ if and
only  if any  subgroup $N$ with $K_{G:A}\neq N\unlhd G$ acts transitively on
$G\!:\!A.$ In particular, a transitive permutation group is quasi-primitive if all its
normal subgroups $\neq 1$ are transitive.

\begin{prop}
\label{1000}
Let $G$ be quasi-primitive on  $G\!:\!A.$  If $A\sm G$ then
$G/\ka$ is simple. Equivalently, if $(G,\Omega)$ is a
quasi-primitive \mt permutation group  then $G$ is simple.
\end{prop}
\pf Suppose that $A\sm G.$ If $G\unrhd N\supset \ka$ then $N$ is
transitive on $G\!:\!A$ as $G$ is quasi-primitive on $G\!:\!A.$
Hence $N=G$ as $A\sm G.$\dne

\vspace{5mm}
\subsection{\sc A Reduction Theorem}

When studying minimal transitivity it is obviously useful to reduce a minimally transitive action $A\sm G$ to one of a smaller group or to one of smaller degree.   Minimal transitivity lends itself to good reduction arguments of this kind for
normal subgroups. For this let $G$ be an arbitary  finite group
with  an \m-stabilizer $A\sm G$ of index $n$ in $G.$ Let $H$ be a
normal subgroup of $G$ with $\ka\subset H$ and $\ka\neq H\neq G.$
Then $H$ is not transitive on $\Omega:=G\!:\!A$ and the  orbits of
$H$ on $\Omega$ are a system of imprimitivity for $G.$ So these
are of the shape $\Omega_{1},...,\Omega_{n^{*}}$ with
$|\Omega_{i}|=s$ and $n^{*}:=\frac{n}{s}.$ Let therefore
$\Omega^{*}:=\{\,\Omega_{i}\,\,|\,\,i=1,...,n^{*}\,\}.$

Let $\Omega_{1}$ be such that it contains the coset $1A$ and let
$B$ be the set-stabilizer of $\Omega_{1}.$ In other words, $B=AH$
and in particular $H\subseteq K_{G:B}.$ Now note  that $B\sm G.$
For if  $M\subseteq G$ satisfies $BM=G$ then $AHM=G.$ As $A\sm G$
we have $G=MH\ka.$ But by choice, $\ka\subseteq H$ so that $G=MH.$
As $H\subseteq K_{G:B}$ we get $G=MK_{G:B}$ and so $B\sm G.$
Equivalently, $G$ acts minimally transitively on $\Om^{*}.$
Therefore we have the following:

\begin{thm} \label{reduct} 
Let $A\sm G$ and suppose $H\neq G$ is normal in $G$ with
$\ka\subset H\neq \ka.$ Then $A\neq AH\sm G.$   \end{thm}

It is worth to formulate this statement in terms of permutation
groups. In conjunction with Proposition~\ref{1000} we have:

\begin{thm} \label{reductperm} 
Let $G$ be a minimally transitive permutation group on $\Omega.$ If $G$ is
quasi-primitive on $\Omega$ then $G$ is simple. Otherwise, if $H$
is a proper normal subgroup of $G$ then  $G$ acts minimally transitively
on the set of $H\!$-orbits. \end{thm}

In other words, a minimally transitive permutation group is either simple or otherwise induces a minimally transitive action on the orbits of any normal subgroup. Another kind of  reduction occurs for the  action of
quotient groups, and this will be used later.
\begin{lem} \label{reduct2} %\marginpar{reduct2}
Let  $N$ be a normal subgroup of  $G$ and let $N\subseteq
A\subseteq G.$ Then $A/N\sm G/N$ if and only if $A\sm G.$
\end{lem}

\pf If  $N\subseteq A\subseteq G$ then $K_{G/N:A/N}=\ka/N.$
Suppose that $A/N\sm G/N$ but  $A\not \sm G.$ So there exist
$H\subseteq G$ with $G=AH$ and $G\neq H\ka.$ Consider
$G/N=A/N\,\cdot\,H\ka/N$ and evaluate $H\ka/N\,\cdot
\,K_{G/N:A/N}=H\ka/N\,\cdot\, \ka/N=H\ka/N\neq G/N,$ a
contradiction. Conversely,   suppose that $A\sm G$ but that
$A/N\not \sm G/N.$ So there is a subgroup $N\subseteq H\subseteq
G$ with $G/N=A/N\,\cdot H/N$ and $H/N\,\cdot \,\ka/N=H\ka/N\neq
G/N.$ So $G=AH$ with $H\ka\neq G,$ a contradiction.\dne

 \vspace{15mm}

\section{\sc Solvable Groups}

For the remainder of the paper we shall restrict ourselves to
minimally transitive representations of solvable groups. If $n$ is
an integer let $\pi(n)$ be the set of primes dividing $n.$ Similarly, $\pi(G)$ and $\pi(G\!:\!H)$ are the prime divisors in  $|G|$ and $|G\!:\!H|$ respectively. Also, $|n|_{p}$ is the highest $p\!$-power dividing $n.$

The following theorem states the basic relation between $\pi(G)$ and the degree of any faithful minimally transitive action when  $G$ is solvable.  For nilpotent groups it
completely characterizes all \mt
actions.

\begin{thm}
\label{1012}
%\marginpar{1012}
\quad (i) Let  $A\sm G$ such that $G/\ka$ is solvable. 
Then $\pi(G\!:\!A)=\pi(G/\ka).$ In particular,
for a solvable minimally transitive permutation group $G$ of
degree $n$ we have $\pi(G)=\pi(n).$

(ii) Let $A\subset G.$ If $A/\ka$ is contained in the Frattini 
subgroup of $G/\ka$ then $A\sm G.$ Conversely, if $G/\ka$ is 
nilpotent and $A\sm G$
then $A/\ka$ is contained in the Frattini subgroup of $G/\ka.$

(iii) If $A\sm G$ and
$|G\!:\!A|=p^{i}$ for some prime $p$ then $G/\ka$ is a $p\!$-group
and $A/\ka$ is contained in the Frattini subgroup of $G/\ka.$
\end{thm}

\pf (i) Let $H$ be a Hall $\pi(G\!:\!A)\!$-subgroup of $G$. Then $AH=G$ as the
left hand side has order $|G|.$ As $A\sm G$ therefore $H\ka=G.$ As
$\ka\subseteq A$ therefore $\pi(G\!:\!A)=\pi(G:\ka).$

(ii) Suppose that  $A\not\sm G.$ Then there exists some maximal subgroup $H'\supseteq \ka$ with $G=AH'$ and $H'\neq G.$ If in addition $A/\ka$ is contained in the Frattini subgroup of $G/\ka$ we have  $A\subseteq H',$ a contradiction. 
Conversely,  if $G/\ka$ is nilpotent and if $H\supseteq \ka$ is  a maximal subgroup of $G$  then $H$ is normal in $G.$ Therefore  $AH\supseteq \ka$ is a group
and  if $A\sm G$ then $AH\neq G,$ and hence 
$A\subseteq H.$ 
(iii) This follows from (i) and (ii). \dne

\medskip

The next result is
a general splitting principle reducing representations of non-nilpotent groups to representations of subgroups, generally involving fewer primes. We denote the  Fitting subgroup of  $X$  by $F(X).$

\begin{thm}
\label{3001}
%\marginpar{3001}
Let $G$ be a solvable group, suppose that $A\sm G$ is core-free and that  $A$ is contained in $F=F(G).$ Let $\pi^{*}:=\pi(G:F)$ and let $Q$ be a Hall $\pi^{*}\!$-subgroup of $G.$ Suppose that $P$ is a normal Sylow $p\!$-subgroup of $G,$ let $A_{P}:=A\cap P$  and $A_{Q}:= A\cap Q.$   Then $p$ does not belong to  $\pi^{*}$. Furthermore,  $A_{Q}\times A_{P}$ is core-free in $Q^{*}P$ and $A_{Q}\times A_{P}\sm Q^{*}P$  for any conjugate $Q^{*}$ of $Q.$

Conversely, let $P_{1},\,P_{2},..,\,P_{t}$  be the normal  Sylow $p_{i}\!$-subgroup of $G.$ Suppose there exist a subgroup $A_{Q}$ of  $F\cap Q$ and subgroups $A_{P_{i}}\subseteq P_{i}$ such that  $ A_{Q}\times A_{P_{i}}\sm Q^{*}P_{i}$ is core-free in $Q^{*}P_{i},$ 
 for all conjugates $Q^{*}$ of $Q$ and all $i=1,..,\,t.$ Then $A_{Q}\times A_{P_{1}}\times \cdots \times A_{P_{t}}\sm G$ is core-free in $G.$
\end{thm}

\medskip For instance, in the simplest case when $\pi^{*}=\{q\}$  we may take $p$ to be any prime in $\pi(n)\setminus \{q\}$ where $n=|G:A|.$ Then $Q\cap F$ is the Sylow $q\!$-subgroup of $F$ so that $A_{Q}$ is the Sylow $q\!$-subgroup of $A.$ Similarly,  $A_{P}$ is the Sylow $p\!$-subgroup of $A$ and hence $A_{Q}\times A_{P}\sm QP$ is minimally transitive of degree $|n|_{q}|n|_{p}.$ Note, for at least one choice of $p$ the group $QP$ is not nilpotent, and evidently, groups of this type are at the basis of any induction in this case. 

\medskip
\pf Evidently, as $P$ is a normal Sylow subgroup of $G$ we have $P\subseteq F.$ Put $K:=K_{QP:(A_{Q}\times A_{P})}.$ Then $K$ is centralized by every Sylow $r\!$-subgroup of $G$ for $r\neq p$ not dividing the order of $Q.$ Further, it is normalized  by $QP$ and hence $K\subseteq A$ is a normal subgroup of $G.$ As $A$ is core-free, $K$ is trivial. 
Now, suppose $A_{Q}\times A_{P}\not\sm QP.$  Then there exists a subgroup $Y\subseteq 
QP$ such that $(A_{Q}\times A_{P})Y=QP$ but 
$Y\neq QP.$ Let $S$ be the direct product of all Sylow $r\!$-subgroups of $F$ for $r\neq p$ not dividing the order of $Q.$ This group is characteristic in $F$ and so normal in $G.$ Therefore, $YS$ is a group and $A(YS)=(AY)S\supseteq(QP)S=QF=G.$ However,  $YS\neq G.$ This is a contradiction, since $A\sm G.$ Finally note that $Q\cap F$ is normal in $F.$ Thus if $Q$ is replaced by $Q^{f}$ then  $A_{Q^{f}}\times A_{P}\sm Q^{f}P.$ But as $A\subseteq F$ we have $A\cap Q=A\cap Q^{f}.$

 Conversely, let $A= A_{Q}\times A_{P_{1}} \cdots \times A_{P_{t}}.$ Since $ A_{Q}\times A_{P_{i}}$ are core-free in $QP_{i}$ also   $A$ is core-free in $G.$ To show that $A\sm G$ suppose that this is not the case. Let therefore $Y$ be a subgroup such that $G=AY$ but $Y\neq G.$ Thus $Y=Y^{*}(Y_{1}\times \cdots \times Y_{t})$ for $Y_{i}:=Y\cap P_{i}$ and  $Y^{*}$  a Hall $\pi^{*}\!$-subgroup  of $Y.$ 
Therefore $Y^{*}\subseteq Q^{f}$ for some  $f\in F.$ 

Further, $G=(P_{1}\times \cdots \times P_{t})Q^{f}$ and since  $Y=(Y_{1}\times \cdots \times Y_{t})Y^{*}$ we have $G=AY=(A_{Q}\times A_{P_{1}} \cdots \times A_{P_{t}})(Y_{1}\times \cdots \times Y_{t})Y^{*}.$ As $A_{Q}$ is a $\pi^{*}\!$-subgroup of $F$ it centralizes all terms other than $Y^{*}.$ Similarly, $Y_{i}$ centralizes all terms $A_{P_{j}}$ with $i\neq j.$ Therefore we can rewrite this as $(P_{1}\times \cdots \times P_{t})Q^{f}= (A_{P_{1}}Y_{1})\times \cdots \times( A_{P_{t}}Y_{t})(A_{Q}Y^{*}).$ For order reasons we have   $A_{Q}Y^{*}=Q^{f}$ and 
$A_{P_{i}}Y_{i}=P_{i}$ for $i=1,..,\,t.$
As $Y \neq G$ we have $Y\cap Q^{f}P_{r}=Y^{*}Y_{r}\neq Q^{f}{P_{r}}$ for at least one $r,$ say $r=1.$ 
Now consider  $(A_{P_{1} }\times A_{Q})(Y^{*}Y_{1})
=(A_{P_{1}}Y_{1})(A_{Q}Y^{*})=P_{1}Q^{f}$
This is a contradiction, since  $A_{P_{1}}\times A_{Q}\sm P_{1}Q^{f}.$ \dne

\vspace{7mm}

{\sc Representations of Square-Free Degree:} \, Next we turn to representations of square-free degree.  Here we get precise information on the Fitting subgroup. 

\begin{thm}\label{square}%\marginpar{square}
For a solvable group  $G$ suppose that $A\sm G$ is core-free and has
square-free index $n$ in $G.$ Let $F$ be the Fitting group of $G.$ Then $|F|$ is coprime  to $|G:F|$ and all Sylow subgroups of $F$ are elementary abelian. In particular, $G$ is nilpotent if and only if $G$ is cyclic of order $n,$ with $A=1.$

Let $\pi^{*}=\pi(n)\setminus \pi(F)$ and let $n^{*}$ be the product of the primes in $\pi^{*}.$ If $C$ is a Hall $\pi^{*}\!$-subgroup of $A$ and if $Q$ is a Hall $\pi^{*}\!$-subgroup of $G$
containing $C$ then $|Q\!:\!C|= n^{*} $ and the action of $Q$ on $Q\!:\!C$ is permutationally equivalent to the action of $G$ on $G\!:\!AF.$

\end{thm}

\medskip
\pf If $n=p_{1}p_{2}\cdots p_{t}$  with pairwise distinct primes
$p_{i}$ then  $\pi(G)=\{p_{1},\,p_{2},\,\dots, p_{t}\}$ by
Theorem~\ref{1012}. Let $N\neq 1$ be a  $p\!$-subgroup of  $G,$
say $p=p_{1},$ which is normal in $G.$ We claim that $N$ is a Sylow subgroup of $G.$ To prove this note that   $N$ has $m:=\frac{n}{p}$ orbits
$\Om_{1},..,\Om_{s},..,\Om_{m}$ on $\Om:=G\!:\!A,$ all of length $|\Om_{s}|=|N\!:\!N\cap A|=p.$

Let $S$ be a Sylow $p\!$-subgroup of $AN.$ As $AN$ is the setwise
stabilizer of the orbit $\Om_{s}$ that contains $1A$ we have that
$p$ does not divide $|G\!:\!AN|.$ Hence  $S$ is a Sylow
$p\!$-subgroup of $G.$  If $Q$ is a Hall $p'\!$-subgroup of $G$ then
$SQ=G$  so that in particular  $G=(AN)Q=A(NQ)$ for order reasons.
As $A\sm G$ we have $(NQ)\ka=G$ but as $A$ is core-free we have
$NQ=G.$ Therefore  $N=S$ is a Sylow subgroup of $G.$ For any $p\in \pi(F)$ let now   $N$ be the unique Sylow $p\!$-subgroup of $F.$ Thus  $N$ is normal in $G$ and hence is a Sylow $p\!$-subgroup of $G.$  It follows that $|F|$ is co-prime to $|G\!:\!F|.$ By the same argument $N$ is characteristically simple and hence elementary abelian. Evidently, if $G=F$ is nilpotent then $G$ is abelian, hence regular on $\Om$ and so cyclic of order $n=|\Om|.$ 

As $|F|$ is co-prime to $|G\!:\!F|$ we may assume for the remainder  that $Q$ is a $\pi^{*}\!$-subgroup of $G$ complementing $F,$ with the further property that  $Q\cap A=C$ is a $\pi^{*}\!$-subgroup of $A.$ Then $G=QF$ with $Q\cap F=1$ and $AF=CF$ with
$C\cap F=1$ implies that the action of $G$ on the $n^{*}$ 
cosets of $AF$ in $G$ is permutationally equivalent to the action of $Q$
on the cosets of $C.$ Hence $C\sm Q$ by Theorem~\ref{reduct}.
 \dne

\bigskip
Some comments are in order. (1) \, While  the
theorem could be formulated for permutation groups the resulting
action of $G$ on the cosets of $AN$ is not faithful, and the same
may be true for the action of $Q$ on the cosets of $C.$\newline
(2) \, As $G$ is solvable there is  at least one normal
$p\!$-subgroup $N,$ as in the proof, and for this $p$ it is the unique normal
$p\!$-subgroup. This subgroup is elementary abelian, and $G$ acts
irreducibly on it. \newline (3) \, The basis of induction for square free degrees occurs when  $n$ is the product of two
distinct primes. A
complete analysis of the possibilities for $G$  can be found in
Suprunenko~\cite{Sup1} and Kopylova~\cite{Kopy1}. For the reader's benefit we collect their results here. 

\begin{thm}\label{sk} 
(Suprunenko~\cite{Sup1}) The permutation group $G$ is minimally transitive of degree $ pq,$ where   $q<p$ are prime numbers with $q$ not dividing
$(p-1),$ if and only if $G$ is isomorphic to one of the following\newline
(i) \,\,\, the cyclic group of order $pq,$\newline
(ii) \,\, a minimal non-abelian
group  G=$PQ,$ where $|Q|=q$ and $P$ is normal in $G,$ with $|P|=p^{m}$ where $m$ is the exponent of  $p$ mod $q,$ or\newline
(iii) \, a minimal non-abelian group   $G=PQ$ with  $Q$ is normal, $|P|=p$ and
$|Q|=q^{r}$ where  $r$  is the exponent of
 $q$ mod $p.$  
\end{thm}

The remaining case where  $q^{r}$ with $r>0$ is the highest power of $q$ dividing  $ (p-1)$ is analyzed in Kopylova~\cite{Kopy1}. Here a similar description is obtained and it is shown  that $G$ is 
 (i) a group of order $q^tp$ with  $0<t\leq r$; \,\, (ii) a group of  order $q^{r+1}p^q$ or \,\, (iii) a group of  order $pq^l$ where $l
$ is the exponent of $q$ mod $p.$

\vspace{7mm}

{\sc Representations of Degree Involving Two Primes:} \,From the discussion so far it is clear that $\{p,\,q\}\!$-groups and their minimally transitive representations play a special role. So let $A\sm G$ be core-free with $\pi(G)=\{p,\,q\}.$  From Theorem~\ref{reduct} it is clear that any normal subgroup $N$ in $G$ gives rise to a minimally transitive representation $AN/N\sm G/N$ of degree $\leq |G:A|.$ Our first observation is the following

\begin{lem}   \,
\label{4009}%\marginpar{4007}
Let  $G$ be solvable and let $A\sm  G$ be core-free in $G.$ Suppose that the prime $q$ divides $|G\!:\!A|$ to the first power only. Let
$N$ be a $q\!$-group which is normal in $G.$ Then $N$ is elementary abelian and is a Sylow subgroup, with $G$ acting irreducibly on $N.$ \end{lem}

\medskip
\pf Let $Q$ be a Sylow
$q\!$-subgroup containing $N$ and let  $P$ be  a $q'\!$-complement in  $G=PQ.$ Note that the $N\!$-orbits on
$G\!:\!A$ are blocks of imprimitivity, all of equal size $q$
and $AN$ is the stabilizer of the orbit containing $1A.$ Therefore $q$ and 
$|G\!:\!AN|$ are co-prime  so that $AN$ contains some Sylow $q\!$-subgroups
of $G.$ From this  we have $(AN)P=G,$ for order reasons. Since
$A\sm G$ and $A(NP)=G$ we have that $(NP)\ka=G$ but $\ka=1$ means
$NP=G.$ This says that $N$ is a Sylow $q\!$-subgroup of $G$ and hence $N=Q.$
Next replace $N$ by a  minimal normal subgroup of $G.$ This group is elementary abelian. By the same argument $N$ has to be Sylow $q\!$-subgroup of $G.$\dne

\bigskip The lemma suggests that the natural choice for a normal subgroup is indeed the Fitting subgroup of $G.$ We follow through with this process when 
$A$ has index $p^{i}q.$ In this case, if either Sylow subgroup $S$ of $G$ is normal then $AS\sm  
G$ gives a minimally transitive representation of a nilpotent group, and this situation is known from Theorem~\ref{1012}. 

Otherwise none of the Sylow subgroups are  normal and by Lemma~\ref{4009} $F_{1}=F(G)$ is a $p\!$-group. By Theorem~\ref{reduct} we have $AF_{1}\sm G$ and if  $K\supseteq F_{1}$ is the core of $AF_{1}$ in $G$ then  $AF_{1}/K\sm G/K$ is minimally transitive, faithful of degree $p^{j}q$ for  $j<i$ and with $|G/K|_{p}<|G|_{p}.$ If  $F_{2}$ is  the pre-image of $F(G/F)$ in $G$ then  $ F_{2}/F_{1}$ is a $q\!$-group. Thus, if $F_{2}$ is not contained in $K$ then Lemma~\ref{4009} shows that $F_{2}K/K$ is a normal Sylow $q\!$-subgroup of $G/K.$  In this case we are reduced to the nilpotent case. 
Otherwise $K\supseteq F_{2}$ so that $|G/K|_{q}<|G|_{q}.$ The process stops when the group becomes  nilpotent or when it is of Suprunenko-Kopylova type.

\end{document}